\documentclass{article}
\usepackage{latexsym,amsfonts,amsmath,amsthm,amssymb,makeidx}
\usepackage[title]{appendix}
\usepackage{CJK,CJKnumb,CJKulem,times,dsfont,ifthen,mathrsfs,latexsym,amsfonts, color}
\usepackage{amsmath,amsthm,makeidx,fontenc,amssymb,bm,graphicx,psfrag,listings, curves,extarrows}
\usepackage{cite}

\let\oldbibliography\thebibliography
\renewcommand{\thebibliography}[1]{%
\oldbibliography{#1}%
\setlength{\itemsep}{0pt}%
}

\makeindex
\newtheorem{definition}{Definition}[section]
\newtheorem{theorem}{Theorem}[section]
\newtheorem{lemma}{Lemma}[section]
\newtheorem{corollary}{Corollary}[section]

\newtheorem{remark}{Remark}[section]

\newcommand{\bt}{\begin{theorem}}
\newcommand{\et}{\end{theorem}}
\newcommand{\bl}{\begin{lemma}}
\newcommand{\el}{\end{lemma}}
\newcommand{\bd}{\begin{definition}}
\newcommand{\ed}{\end{definition}}
\newcommand{\bc}{\begin{corollary}}
\newcommand{\ec}{\end{corollary}}
\newcommand{\bp}{\begin{proof}}
\newcommand{\ep}{\end{proof}}
\newcommand{\bx}{\begin{example}}
\newcommand{\ex}{\end{example}}
\newcommand{\bi}{\begin{exercise}}
\newcommand{\ei}{\end{exercise}}
\newcommand{\bo}{\begin{prop}}
\newcommand{\eo}{\end{prop}}
\newcommand{\br}{\begin{remark}}
\newcommand{\er}{\end{remark}}
\newcommand{\be}{\begin{equation}}
\newcommand{\ee}{\end{equation}}
\newcommand{\ba}{\begin{align}}
\newcommand{\ea}{\end{align}}
\newcommand{\bn}{\begin{enumerate}}
\newcommand{\en}{\end{enumerate}}
\newcommand{\bg}{\begin{align*}}
\newcommand{\bcs}{\begin{cases}}
\newcommand{\ecs}{\end{cases}}

\newcommand{\bean}{\begin{eqnarray*}}
\newcommand{\eean}{\end{eqnarray*}}

\numberwithin{equation}{section}

\begin{document}
\title{{\bf  Asymptotic behavior of  positive solutions  to a nonlinear biharmonic equation near isolated singularities}}

\date{}
\author{\\{\bf Hui  Yang}\thanks {E-mail addresses:  mahuiyang@ust.hk; hui-yang15@outlook.com }
\\
\footnotesize {\it  Yau Mathematical Sciences Center, Tsinghua University, Beijing 100084, China}\\
}

\maketitle
\begin{center}
\begin{minipage}{120mm}
\begin{center}{\bf Abstract}\end{center}

In this paper, we consider the asymptotic behavior of  positive solutions of  the biharmonic equation  
$$ \Delta^2 u = u^p~~~~~~~\textmd{in} ~ B_1\backslash \{0\}$$   
with an isolated singularity, where the punctured  ball $B_1 \backslash \{0\} \subset \mathbb{R}^n$ with $n\geq 5$ and $\frac{n}{n-4}  <  p  < \frac{n+4}{n-4}$. This equation is relevant for the $Q$-curvature problem in conformal geometry. We classify isolated singularities of positive solutions and describe the asymptotic behavior of positive singular solutions  without the sign assumption for $-\Delta u$.    We also give a new method to prove removable singularity theorem for nonlinear higher order equations.  


\vskip0.10in


\end{minipage}

\end{center}

\vskip0.390in

\section{Introduction and main results}
In this paper, we study the asymptotic behavior of  positive solutions of  the biharmonic equation  
\begin{equation}\label{Bi}
 \Delta^2 u = u^p~~~~~~~\textmd{in} ~ B_1\backslash \{0\}
\end{equation}
with an isolated singularity, where the punctured  ball $B_1 \backslash \{0\} \subset \mathbb{R}^n$ with $n\geq 5$ and $\frac{n}{n-4}  <  p  < \frac{n+4}{n-4}$. Here the unit ball $B_1$ can be replaced by any bounded domain $\Omega \subset \mathbb{R}^n$ containing 0.   This equation serves as a basic model of nonlinear fourth-order equations and is also related to the $Q$-curvature problem in conformal geometry.  Equation \eqref{Bi}  and related equations  arise in several models describing various phenomena in the applied sciences see, for instance, Gazzola, Grunau and Sweers \cite{GGS}.   For an introduction to the $Q$-curvature problem see, for instance,  Hang and Yang \cite{HY}.   

\vskip0.10in

We first recall that the corresponding second order equation (when $n\geq 3$ and $\frac{n}{n-2} < p < \frac{n+2}{n-2}$) 
\begin{equation}\label{Sec}
-\Delta u = u^p~~~~~~~\textmd{in} ~ B_1\backslash \{0\}
\end{equation}
was studied by Gidas-Spruck \cite{GS} and Caffarelli-Gidas-Spruck \cite{CGS}. More specifically, the following classification result is obtained.\\

\noindent{\bf Theorem A}  (\cite{CGS,GS}) Let $n \geq 3$ and $u \in C^2(B_1 \backslash \{0\})$  be a positive solution of \eqref{Sec}. Assume 
$$\frac{n}{n-2} < p < \frac{n+2}{n-2}.$$
Then either the singularity at $x=0$ is removable, or $u$ is a distribution solution in the entire ball $B_1$, and 
$$
\lim_{|x| \to 0} |x|^{\frac{2}{p-1}}u(x)=C_0,
$$
where
$$
C_0= \left\{ \frac{2(n-2)}{(p-1)^2} \left( p-\frac{n}{n-2} \right) \right\}^{\frac{1}{p-1}}. 
$$
\\
In addition, the asymptotic behavior of positive solutions of \eqref{Sec} near an isolated singularity was studied by Lions \cite{L}  for $1 < p < \frac{n}{n-2}$,  by Aviles \cite{A} for $p=\frac{n}{n-2}$, by Caffarelli-Gidas-Spruck\cite{CGS} and by Korevaar-Mazzeo-Pacard-Schoen \cite{KMPS} in the case $p=\frac{n+2}{n-2}$ and by Bidaut-V\'{e}ron and  V\'{e}ron \cite{BV} when $p>\frac{n+2}{n-2}$. Hence the isolated singularities of positive solutions for the second order equation \eqref{Sec} have been very well understood.  The asymptotic behavior of positive solutions for a more general second order equation $-\Delta u =K(x) u^{\frac{n+2}{n-2}}$ with isolated singularity was studied by Chen-Lin \cite{CL1,CL2} and Taliaferro-Zhang \cite{TZ}. See also Gonz\'alez \cite{Gon}, Li \cite{Li} and Han-Li-Teixeira \cite{HLT}  for a fully nonlinear equation of second order. 

\vskip0.10in

In the fundamental paper \cite{Lin}, Lin classified all positive smooth entire solutions of \eqref{Bi} with $1< p \leq \frac{n+4}{n-4}$ in $\mathbb{R}^n$ via the moving plane method. We refer to Chang-Yang\cite{CY}, Martinazzi \cite{Ma} and  Wei-Xu \cite{WX} for the  classification of smooth solutions of  the higher-order equations in $\mathbb{R}^n$.    For the supercritical case, that is for $p>\frac{n+4}{n-4}$, the positive smooth radial symmetric  solutions of \eqref{Bi} in $\mathbb{R}^n$ were studied by Gazzola-Grunau \cite{GG}, Guo-Wei \cite{GW} and Winkler \cite{W}.  We also refer to a recent paper Frank-K\"onig \cite{FK} for a classification of positive singular solutions to \eqref{Bi} with $p=\frac{n+4}{n-4}$ in $\mathbb{R}^n \backslash \{0\}$, where the positive singular solutions are radially symmetric about the origin (see Theorem 4.2 in \cite{Lin}).    

\vskip0.10in

As far as we know, the classification of isolated singularities of positive solutions and the asymptotic behavior of positive singular solutions to fourth order equation \eqref{Bi} in $B_1\backslash \{0\}$ are far less known than the second order problem \eqref{Sec}. Remark that, a positive solution $u$ of \eqref{Bi} in $B_1\backslash \{0\}$  may not be radially symmetric.

\vskip0.10in

If one looks closely at the tools being used in the proofs of second order  problems, then one finds that the maximum principle plays an essential role. This is a crucial distinction from higher order problems for which there is no the maximum principle.  Here and in the sequel "higher order" means order at least four. Another important tool intensively used for second order problems is the truncation method. This method is powerful in regularity theory and in properties of first order Sobolev spaces. However, the truncation method also fails for higher order problems.  Therefore, the methods of  above mentioned papers for second order problems cannot be applied to the fourth order equation \eqref{Bi}.

\vskip0.10in

Nevertheless we succeed here in proving exact asymptotic behavior of positive singular solutions for \eqref{Bi} which is completely analogous to its second order counterpart.  Remark that our proof is very different from that of Theorem A in \cite{CGS,GS}.   
Our main result is the following

\begin{theorem}\label{T-Bi}
Let $n \geq 5$ and $u \in C^4(B_1 \backslash \{0\})$  be a nonnegative  solution of \eqref{Bi}. Assume 
$$\frac{n}{n-4} < p < \frac{n+4}{n-4}.$$
Then either the singularity at $x=0$ is removable, or $u$ is a distribution solution in the entire ball $B_1$, and 
\begin{equation}\label{Beh}
\lim_{|x| \to 0} |x|^{\frac{4}{p-1}}u(x)=C_{p, n} >0,
\end{equation}
where $ C_{p,n}= \left[ K_0(p, n) \right]^{\frac{1}{p-1}}$ and 
\begin{equation}\label{K0}
\aligned
K_0(p, n)= & \frac{8}{(p-1)^4}\Big[ (n-2)(n-4)(p-1)^3 + 2(n^2 -10n +20)(p-1)^2 \\
& -16(n-4)(p-1) + 32 \Big].  
\endaligned
\end{equation}
\end{theorem}

\begin{remark}
We don't need any additional assumptions for $-\Delta u$ in $B_1\backslash \{0\}$ and for boundary conditions.  Soranzo \cite{Sor} studied the local behavior of positive solutions of \eqref{Bi} with additional assumption 
\begin{equation}\label{Lap}
-\Delta u \geq 0   ~~~~~~~ \textmd{in} ~  B_1\backslash \{0\}.
\end{equation} 
Under the assumption \eqref{Lap},  Soranzo \cite{Sor}  classified the isolated singularities of positive solutions of \eqref{Bi} for $1 < p <\frac{n}{n-4}$ and obtained an upper bound of  radially  symmetric positive solutions of \eqref{Bi} for $p \geq \frac{n}{n-4}$. Theorem \ref{T-Bi} also answers an open question raised in \cite{Sor} (see Remark 5 there) and shows, in particular, that the nonnegativity of $-\Delta u$ in this problem is not necessary. 
Recently, Jin and Xiong \cite{JX} proved sharp blow up rates and the asymptotic radial symmetry of  positive solutions of \eqref{Bi} with $p=\frac{n+4}{n-4}$ near the singularity  under the sign assumption \eqref{Lap}.   We also mention that Ferrero-Grunau \cite{FG} have obtained the asymptotic behavior of positive radial singular solutions for biharmonic operator and power-like nonlinearity with the Dirichlet boundary condition.  
\end{remark}

\begin{remark}
In $\mathbb{R}^n$, suppose $u$ is a positive smooth function satisfies equation \eqref{Bi} with $p>1$, then necessarily we have
 $$
 -\Delta u > 0 ~~~~~~ \textmd{in} ~ \mathbb{R}^n. 
 $$
See Theorem 3.1 in   Wei-Xu \cite{WX}. This important fact about $-\Delta u$ enables the maximum principle to be applied to positive solutions  of  \eqref{Bi} in $\mathbb{R}^n$.  Such as see \cite{Lin,WX}.  Hence the positive solutions of equation \eqref{Bi} in $\mathbb{R}^n$  provide  enough information for applying the maximum principle, but this is not true for \eqref{Bi} in $B_1 \backslash \{0\}$.  
\end{remark}

\begin{remark}
When $\frac{n}{n-4} < p < \frac{n+4}{n-4}$, it is well known that the function 
$$u(x)= C_{p,n}|x|^{-\frac{4}{p-1}}$$
is an exact positive singular solution of \eqref{Bi} which obviously satisfies asymptotic behavior \eqref{Beh}.  See also Guo-Wei-Zhou \cite{GWZ} for another a family of positive singular radial solutions  of \eqref{Bi} in $\mathbb{R}^n \backslash \{0\}$. 
\end{remark}

The rest of this paper is organized as follows. In Section 2, we establish some basic estimates.  In Section 3, we prove Theorem \ref{T-Bi}.

\section{Preliminaries}

In this section we establish some basic estimates.  First we recall the following Liouville type theorem.  For its  proof, such as see Lin \cite{Lin}. 

\begin{theorem}\label{T201}
(\cite{Lin}) Suppose that $u$ is a nonnegative solution of
\begin{equation}\label{Lio}
\Delta^2 u=u^p~~~~~~~ \textmd{in} ~ \mathbb{R}^n
\end{equation}
for $1 < p < \frac{n+4}{n-4}$. Then $u \equiv 0$ in $\mathbb{R}^n$. 
\end{theorem}
By a doubling lemma of  Pol\'acik, Quittner and Souplet \cite{PQS} and above Liouville theorem, we have the following singularity and decay estimates. Because their proof is similar, we only give the proof of decay estimates here. 
\begin{lemma}\label{Le201}
Let $u \in C^4(B_1 \backslash \{0\})$  be a nonnegative solution of \eqref{Bi} with $1< p< \frac{n+4}{n-4}$. Then 
\begin{equation}\label{L2201}
u(x) \leq C |x|^{-\frac{4}{p-1}} ~~~~~~ for ~|x| \leq \frac{1}{2},
\end{equation}
where $C$ is a constant, depending on $n$ and $p$ only. 
\end{lemma}

\begin{remark}
For the second order equation \eqref{Sec}, if one has an upper estimate similar to \eqref{L2201}, then one can easily obtain the following  Harnack inequality
\begin{equation}\label{L2401}
\sup_{r \leq |x| \leq 2r} u  \leq C \inf_{r \leq |x| \leq 2r} u, 
\end{equation} 
where $C$ is independent of $r$. Such as see \cite{A,CL2,KMPS,GS}. This  is an essential tool  for these papers to study isolated singularities of second order problems.  In a recent paper \cite{CJSX} Caffarelli, Jin, Sire and Xiong  use a similar Harnack inequality to classify isolated singularities of positive solutions of a fractional equation.  However, this Harnack inequality does not  generally hold for fourth order equation \eqref{Bi}.  In particular, if we suppose additionally that  
\begin{equation}\label{Nen}
-\Delta u \geq 0 ~~~~~~~\textmd{in} ~ B_1\backslash \{0\},
\end{equation}
then Caristi-Mitidieri \cite{CM}  proved that the similar Harnack inequality still  holds for fourth order equation \eqref{Bi}. 
\end{remark}

\begin{remark}
We also remark that the condition \eqref{Nen} is necessary for the validity of the Harnack inequality to  biharmonic equations as the following simple example shows: consider the function $u(x)=\sum_{i=1}^n x_i^2$. It is nonnegative, satisfies $\Delta^2 u=0$ and $\Delta u =2n$, but the Harnack inequality does not hold in $B_1(0)$. 
\end{remark}

\begin{lemma}\label{Le2089}
Let $u$ be a nonnegative solution of 
\begin{equation}\label{Inf89}
\Delta^2 u = u^p~~~~~ \textmd{in} ~ B_1^c, 
\end{equation}
where $B_1^c:=\{x \in \mathbb{R}^n : |x| > 1\}$.  Assume $1< p< \frac{n+4}{n-4}$. Then 
\begin{equation}\label{L227801}
u(x) \leq C |x|^{-\frac{4}{p-1}} ~~~~~~ for ~|x| > 2,
\end{equation}
where $C$ is a constant, depending on $n$ and $p$ only. 
\end{lemma}

\bp
Suppose by contradiction that there exist a sequence of nonnegative solutions $(u_k)_k$ of \eqref{Inf89} and a sequence of points $|x_k| >2$,  such that
$$
M_k(x_k) d(x_k) > 2k,~~~~~~~~ k=1, 2, \cdots, 
$$
where $M_k(x) := \left( u_k(x) \right)^{\frac{p-1}{4}}$  and $d(x) := \textmd{dist}(x, \partial B_1^c) =|x| -1$ for $x \in B_1^c$. By the doubling lemma of \cite{PQS} there exists another sequence $y_k \in B_1^c$ such that
$$
M_k(y_k)d(y_k) > 2k,  ~~~~~~ M_k(y_k) \geq M_k(x_k)
$$
and 
$$
M_k(z) \leq 2M_k(y_k) ~~~~\textmd{for} ~ \textmd{any} ~ |z - y_k|\leq k \lambda_k. 
$$
where $\lambda_k:=M_k(y_k)^{-1}$. We now define
$$
\bar{u}_k(x)=\lambda_k^{\frac{4}{p-1}} u_k(y_k + \lambda_k x)~~~~~~ \textmd{for} ~ x\in B_k(0). 
$$
Then $\bar{u}_k$ is a nonnegative solution of
$$
\Delta^2 \bar{u}_k =(\bar{u}_k)^p~~~~~~~~ \textmd{in} ~ B_k(0). 
$$
Moreover,
\begin{equation}\label{M941}
\bar{u}_k(0)=1~~~~~~ \textmd{and} ~~~~~~ \max_{B_k(0)} |\bar{u}_k| \leq 2^{\frac{4}{p-1}}. 
\end{equation}
By the elliptic estimates, we deduce that a subsequence of $(\bar{u}_k)_k$ converges in $C_{loc}^4(\mathbb{R}^n)$ to a nonnegative solution $u_\infty$ of \eqref{Lio} in $\mathbb{R}^n$.  By \eqref{M941}, we have $u_\infty(0)=1$. This contradicts Theorem \ref{T201}. 
\ep

\begin{corollary}\label{Co201}
Let $u \in C^4(B_1 \backslash \{0\})$ be a nonnegative solution of \eqref{Bi} with $1< p< \frac{n+4}{n-4}$. Then there exists a constant $C_1=C_1(n, p)$ such that for all $|x|\leq \frac{1}{4}$,
\begin{equation}\label{L2301}
\sum_{k \leq 3} |x|^{\frac{4}{p-1} +k} |\nabla^k u (x)| \leq C_1. 
\end{equation}
\end{corollary}

\bp
For any $x_0$ with $|x_0|\leq \frac{1}{4}$, take $\lambda=\frac{|x_0|}{2}$ and define
$$
\bar{u}(x)=\lambda^{\frac{4}{p-1}} u(x_0 + \lambda x). 
$$
Then $\bar{u}$ is a nonnegative solution of  \eqref{Bi} in $B_1$. By the Lemma \ref{Le201}, $|\bar{u}| \leq C_2$ in $B_1$. The standard elliptic estimates give
$$
\sum_{k \leq 3} |\nabla^k \bar{u}(0)| \leq C_3.
$$
Rescaling back we obtain \eqref{L2301}. 
\ep

Using a similar scaling argument as above, we also have 

\begin{corollary}\label{Co2089}
Let $u$ be a nonnegative solution of \eqref{Inf89} with $1< p< \frac{n+4}{n-4}$. Then there exists a constant $C_2=C_2(n, p)$ such that for all $|x|\geq 4$,
\begin{equation}\label{L2301}
\sum_{k \leq 3} |x|^{\frac{4}{p-1} +k} |\nabla^k u (x)| \leq C_2. 
\end{equation}
\end{corollary}

\section{Proof of the main result} 
In this section we will prove Theorem \ref{T-Bi}. We first show that any nonnegative solution of \eqref{Bi} with $p \geq \frac{n}{n-4}$  is a solution in $B_1$ in the sense of distribution.

\begin{lemma}\label{Weak}
Assume $p \geq \frac{n}{n-4}$ and  that $u \in C^4(B_1 \backslash \{0\})$  is a nonnegative solution  of \eqref{Bi}. Then $u\in L_{loc}^p(B_1)$ and $u$ is a distribution solution of \eqref{Bi} in $B_1$, that is,
\begin{equation}\label{Dis}
\int_{B_1} u \Delta^2 \varphi = \int_{B_1} u^p \varphi ~~~~~ for ~ all ~ \varphi \in C_c^\infty (B_1). 
\end{equation}
\end{lemma}
\bp
For any  $ 0< \epsilon \ll 1$, we  take $\eta_\epsilon \in C^\infty (\mathbb{R}^n)$ with values in $[0, 1]$ satisfying 
\begin{equation}\label{Eta90}
\eta_\epsilon(x)=
\begin{cases}
0~~~~~~&\textmd{for}~|x|\leq \epsilon,\\
1~~~~~~&\textmd{for}~|x|\geq 2\epsilon
\end{cases}
\end{equation}
and
\begin{equation}\label{Eta91}
|\nabla^k \eta_\epsilon (x)| \leq C \epsilon^{-k} ~~~~~~~~~~~ \textmd{for} ~  k=1, 2, 3, 4. 
\end{equation}
Let $m=\frac{4p}{p-1}$ and define $\xi_\epsilon=(\eta_\epsilon)^m$.   Multiplying \eqref{Bi} by $\xi_\epsilon $ and integrating by parts in $B_r$ with $\frac{1}{2} < r<1$, we get
$$
\int_{B_r} u^p \xi_\epsilon = \int_{\partial B_r}  \frac{\partial}{\partial \nu} \Delta u + \int_{B_r} u  \Delta^2\xi_\epsilon. 
$$
Note that
$$
| \Delta^2\xi_\epsilon |  \leq C \epsilon^{-4}  (\eta_\epsilon)^{m-4}\chi_{\{\epsilon \leq |x|\leq 2\epsilon\}} =C \epsilon^{-4} (\xi_\epsilon)^{1/p} \chi_{\{\epsilon\leq |x|\leq 2\epsilon\}}.
$$
By H\"older's inequality, we obtain
$$
\aligned
\left| \int_{B_r} u  \Delta^2\xi_\epsilon \right| & \leq C \epsilon^{-4} \int_{\{\epsilon\leq |x|\leq 2\epsilon\}} u (\xi_\epsilon)^{1/p} \\
& \leq C \epsilon^{-4} \cdot \epsilon^{n(1-1/p)} \left( \int_{\{\epsilon \leq |x|\leq 2\epsilon\}} u^p \xi_\epsilon \right)^{1/p} \\
& \leq C \left( \int_{\{\epsilon \leq |x|\leq 2\epsilon\}} u^p \xi_\epsilon \right)^{1/p}.
\endaligned
$$
Hence we have
$$
\int_{B_r} u^p \xi_\epsilon \leq  \int_{\partial B_r}  \frac{\partial}{\partial \nu} \Delta u + C \left( \int_{\{\epsilon \leq |x|\leq 2\epsilon\}} u^p \xi_\epsilon \right)^{1/p}.
$$
This implies that there exists a constant $C>0$ (independent of $\epsilon$) such that
$$
\int_{B_r} u^p \xi_\epsilon \leq C.
$$
Now letting $\epsilon \to 0$, we conclude that $u\in L^p(B_r)$.  

\vskip0.10in

To show that $u$ is a distribution solution we need to  establish \eqref{Dis}. For any $\varphi \in C_c^\infty (B_1)$, using $\eta_\epsilon \varphi$ as a test function in \eqref{Bi} with $\eta_\epsilon$ as before gives 
\begin{equation}\label{Dis01}
\int_{B_1} u \Delta^2(\eta_\epsilon\varphi) = \int_{B_1} u^p \eta_\epsilon\varphi. 
\end{equation}
By a direct computation, we have 
$$
\aligned
\Delta^2(\eta_\epsilon\varphi) & =  \eta_\epsilon \Delta^2 \varphi + 4 \nabla \eta_\epsilon \cdot \nabla \Delta\varphi + 2\Delta \eta_\epsilon \Delta\varphi +4 \sum_{i,j=1}^n (\eta_\epsilon)_{x_ix_j} \varphi_{x_i x_i}\\
& ~~~~~ + 4 \nabla \Delta\eta_\epsilon \cdot \nabla \varphi + \varphi \Delta^2 \eta_\epsilon   \\
& =: \eta_\epsilon \Delta^2 \varphi  + \psi, 
\endaligned
$$
and by H\"older's inequality, we get 
$$
\aligned
\left| \int_{B_1} u \psi \right| & \leq C \epsilon^{-4} \int_{\{\epsilon \leq |x| \leq 2\epsilon\}} u \\
& \leq C \epsilon^{-4} \cdot \epsilon^{n(1-1/p)}  \left( \int_{\{\epsilon \leq |x| \leq 2\epsilon\}} u^p \right)^{1/p} \\
& \leq C \left( \int_{\{\epsilon \leq |x| \leq 2\epsilon\}} u^p \right)^{1/p} \to 0 ~~~~~ \textmd{as} ~ \epsilon \to 0. 
\endaligned
$$
Letting $\epsilon \to 0$ in \eqref{Dis01}, then \eqref{Dis} follows immediately from the dominated convergence theorem and the proof is complete. 
\ep

Now we prove that if $u$ is a nonnegative solution of \eqref{Bi}  in $\mathbb{R}^n \backslash \{0\}$, then the sign condition 
 $$
 - \Delta u \geq 0 ~~~~~~~ \textmd{in} ~ \mathbb{R}^n \backslash \{0\}
 $$
holds.  This allows us to use the maximum principle for $u$ in $\mathbb{R}^n \backslash \{0\}$.  

\begin{lemma}\label{Non589}
Assume $\frac{n}{n-4} < p <\frac{n+4}{n-4}$ and  that $u \in C^4(\mathbb{R}^n \backslash \{0\})$ is a nonnegative solution  of 
\begin{equation}\label{WS69}
\Delta^2 u =u^p~~~~~~~ \textmd{in} ~ \mathbb{R}^n \backslash \{0\}. 
\end{equation}
Then $-\Delta u$ is a superharmonic function in  $\mathbb{R}^n$ in the distributional sense. Moreover,
$$
-\Delta u \geq 0 ~~~~~~~ \textmd{in} ~ \mathbb{R}^n \backslash \{0\}. 
$$
\end{lemma}

\bp
By Lemma \ref{Weak}, we have $u \in L_{loc}^p(\mathbb{R}^n)$. Let  $\varphi \in C_c^\infty (\mathbb{R}^n)$ be a nonnegative function. We will  prove that 
$$
\int_{\mathbb{R}^n} \Delta u \Delta \varphi \geq 0. 
$$
Let $\eta_\epsilon \in C^\infty(\mathbb{R}^n)$ satisfy \eqref{Eta90} and \eqref{Eta91}. Multiplying \eqref{WS69} by $\eta_\epsilon \varphi$ and integrating by parts, we obtain 
$$
\aligned
0 & \leq \int_{\mathbb{R}^n} \eta_\epsilon\varphi u^p\\
& = \int_{\mathbb{R}^n}  \Delta(\eta_\epsilon\varphi) \Delta u\\
& = \int_{\mathbb{R}^n} \Delta u ( \Delta \varphi \eta_\epsilon + 2\nabla\varphi \cdot \nabla\eta_\epsilon +\varphi \Delta\eta_\epsilon). 
\endaligned
$$
Denote $\psi=2\nabla\varphi \cdot \nabla\eta_\epsilon +\varphi \Delta\eta_\epsilon$. Then $\psi(x) \equiv 0$ for $|x| \leq \epsilon$ and for $|x| \geq 2\epsilon$, and 
$$
|\Delta \psi (x)| \leq C \epsilon^{-4}.
$$ 
Since $n-4-\frac{n}{p} >0$, we have
$$
\aligned
\left| \int_{\mathbb{R}^n} \Delta u \psi \right| & \leq \int_{\mathbb{R}^n} u |\Delta \psi| \\
& \leq C \epsilon^{-4} \left( \int_{\{\epsilon \leq |x| \leq 2\epsilon\}} u^p\right)^{1/p} \epsilon^{n(1-1/p)} \\
& \leq C \epsilon^{n-4-\frac{n}{p}} \to 0,~~~~~~ \textmd{as} ~ \epsilon \to 0. 
\endaligned
$$
Therefore, we obtain 
$$
\aligned
\int_{\mathbb{R}^n} \Delta u \Delta \varphi & =\lim_{\epsilon \to 0} \int_{\mathbb{R}^n} \Delta u ( \Delta \varphi \eta_\epsilon + 2\nabla\varphi \cdot \nabla\eta_\epsilon +\varphi \Delta\eta_\epsilon) \\
& =  \int_{\mathbb{R}^n} \varphi u^p \geq 0. 
\endaligned
$$
Thus, $-\Delta u$ is a superharmonic  function in  $\mathbb{R}^n$ in the distributional sense.  

\vskip0.10in

Let  $v_\epsilon := -\Delta u +\epsilon$ for $\epsilon >0$. By Corollary \ref{Co2089},  we have $\lim_{|x| \to \infty} |\Delta u(x)| = 0$. Therefore,  for any $\epsilon >0$,  there exists $R_\epsilon$ such that
$$
v_\epsilon > \frac{\epsilon}{2} ~~~~~~~~  \textmd{for} ~ |x|\geq R_\epsilon.
$$
Since $v_\epsilon$ is also a superharmonic function in  $\mathbb{R}^n$ in the distributional sense, we obtain 
$$
v_\epsilon \geq 0~~~~~~~ \textmd{in} ~ \mathbb{R}^n \backslash \{0\}. 
$$
Letting $\epsilon \to 0$, we get $-\Delta u \geq 0$ in $\mathbb{R}^n \backslash \{0\}$. This completes the proof. 
\ep

\vskip0.10in

Let $u$ be a nonnegative solution of \eqref{Bi}. We use the following transformation of \eqref{Bi} (also known as {\it Emden-Fowler transformation}): set
\begin{equation}\label{Set1}
t=\ln |x|,~~~~~ \theta=\frac{x}{|x|}
\end{equation}
and 
\begin{equation}\label{Set2}
w(t, \theta)=|x|^{\frac{4}{p-1}}u(|x|, \theta)=e^{4t / (p-1)} u(e^t, \theta), ~~~~~~ t \in (-\infty, 0), ~ \theta \in \mathbb{S}^{n-1}. 
\end{equation}
By a tedious computation  we find that equation \eqref{Bi} for $u$ is equivalent to the following equation for $w$:
\begin{equation}\label{W01}
\aligned
 \partial_t^{(4)}w &  + K_3 \partial_t^{(3)} w + K_2 \partial_t^{(2)} w + K_1 \partial_t  w +\Delta_\theta^2 w + 2\partial_t^{(2)} \Delta_\theta w \\
&  + K_3 \partial_t\Delta_\theta w + J_1\Delta_\theta w +K_0 w = w^p  ~~~~~ \textmd{in} ~  (-\infty, 0) \times \mathbb{S}^{n-1}, 
\endaligned
\end{equation}
where $\Delta_\theta$ is the Beltrami-Laplace operator on $\mathbb{S}^{n-1}$, the constants $K_i=K_i(p, n)$ $(i=0, \cdots, 3)$ and $J_1=J_1(p, n)$ are given by 
$$
\aligned
K_0  & = \frac{8}{(p-1)^4} \Big[ (n-2)(n-4)(p-1)^3 + 2(n^2 -10n + 20)(p-1)^2 \\
& ~~~~~ -16(n-4)(p-1) +32 \Big], \\
K_1 & = -\frac{2}{(p-1)^3} \Big[ (n-2)(n-4)(p-1)^3 + 4(n^2 -10n +20)(p-1)^2 \\
& ~~~~~ -48(n-4)(p-1) +128 \Big], \\
K_2 & = \frac{1}{(p-1)^2} \Big[ (n^2 - 10n +20)(p-1)^2 -24(n-4)(p-1) +96 \Big], \\
K_3 & = \frac{2}{ p-1 }\Big[ (n-4)(p-1) - 8 \Big], \\
J_1 &= -\frac{2}{(p-1)^2} \Big[ (n-4)(p-1)^2 + 4(n-4) (p-1) - 16 \Big]. \\
\endaligned
$$
Note that  if $ p < \frac{n+4}{n-4}$, then
\begin{equation}\label{Equ}
(n-4)(p-1) < 8. 
\end{equation}
It is not difficult to show that 
\begin{equation}\label{Equ01}
K_1=K_3=0  ~~~~~~~ \textmd{if}  ~~ p=\frac{n+4}{n-4}. 
\end{equation} 
Moreover, we have
\begin{lemma}\label{L-S301}
Assume $n \geq5$ and $\frac{n}{n-4} < p < \frac{n+4}{n-4}$. Then
\begin{equation}\label{Equ02}
K_0 >0, ~~~~~~~~ K_1 >0, ~~~~~~~~ K_3<0.
\end{equation} 
\end{lemma}
\begin{remark}
We emphasize that the sign of $K_1$ and $K_3$ will be essentially used in our arguments. We also point that $J_1 < 0$ for $\frac{n}{n-4} < p <\frac{n+4}{n-4}$ and the sign of $K_2$ depends on $p$ and $n$. 
\end{remark}
\bp
By \eqref{Equ}, we easily obtain $K_3 < 0$. Next we will prove that $K_1 >0$ under the assumptions.  For this purpose, we consider the function
$$
f(s) = (n-2)(n-4)s^3 + 4(n^2 -10n +20)s^2 -48(n-4)s +128 
$$
with $s \in (\frac{4}{n-4}, \frac{8}{n-4})$.  Then
$$
f^\prime(s)= 3(n-2)(n-4) s^2 + 8(n^2 -10n +20) s -48(n-4). 
$$
Since $f^\prime(0) <0$, $f^\prime$ has only one positive root, we denote it by $s_+$. We also denote 
$$
s_0=\frac{4}{n-4}~~~~~~ \textmd{and} ~~~~~~ s_1=\frac{8}{n-4}.
$$
By a direct calculation, we have $f^\prime(s_1)=\frac{16(n^2 - 4n +8)}{n-4} >0$. Hence we must have $s_+ < s_1$. We consider  separately the case $s_0 \geq s_+$ and the case $s_0 < s_+$. 

\vskip0.10in 

{\it Case 1: $s_0 \geq s_+$. } In this case we have $f^\prime(s) > 0$ for all $s\in (s_0, s_1)$. By \eqref{Equ01}, 
$$
f(s) <  f(s_1)=0 ~~~~~~ \textmd{for} ~ \textmd{any} ~ s\in(s_0, s_1). 
$$

{\it Case 2: $s_0 < s_+$.} In this case we have $f^\prime(s) < 0$ in $(s_0, s_+)$ and $f^\prime(s)>0$ in $(s_+, s_1)$. Combining \eqref{Equ01} and the basic fact $f(s_0)=-\frac{64(n-2)}{(n-4)^2} <0$, we obtain 
$$
f(s) < \max\{ f(s_0), f(s_1)\}=0 ~~~~~~ \textmd{for} ~ \textmd{any} ~ s\in(s_0, s_1).  
$$
From these we easily get $K_1 >0$ if $\frac{n}{n-4} < p< \frac{n+4}{n-4}$.   

\vskip0.10in 

Now we check $K_0 >0$. Similarly, we consider
$$
g(s)=(n-2)(n-4)s^3 + 2(n^2 -10n +20)s^2 -16(n-4)s +32
$$
with $s > \frac{4}{n-4}$.  Then
$$
g^\prime(s)= 3(n-2)(n-4) s^2 + 4 (n^2 -10n +20) s - 16(n-4). 
$$
Direct calculations show that $g^\prime(s_0)=\frac{16(n-2)}{n-4} >0$ and $g(s_0)=0$.  From this we get  $g^\prime(s) >0$ for all $s > s_0$ and then
$$
g(s) > g(s_0)= 0  ~~~~~~ \textmd{for} ~ \textmd{any} ~ s > s_0. 
$$
Hence we have $K_0 >0$ if $p>\frac{n}{n-4}$.  
\ep

Next we will establish an important monotonicity formula. Let $w$ be a nonnegative solution of \eqref{W01}. Define
$$
\aligned
E(t; w) := & \int_{\mathbb{S}^{n-1}} \partial_t^{(3)}w \partial_t w  -\frac{1}{2}\int_{\mathbb{S}^{n-1}} \Big[ \left( \partial_t^2 w \right)^2 - 2K_3 \partial_t^2 w \partial_t w - K_2 \left( \partial_t w \right)^2  \Big] \\
& +\frac{1}{2} \int_{\mathbb{S}^{n-1}}\Big[ |\Delta_\theta w|^2 - J_1 |\nabla_\theta w|^2 \Big]  + \frac{K_0}{2} \int_{\mathbb{S}^{n-1}} w^2 \\
& -\frac{1}{p+1}\int_{\mathbb{S}^{n-1}} w^{p+1} - \int_{\mathbb{S}^{n-1}} |\partial_t \nabla_\theta w|^2. 
\endaligned
$$
Then we have the following 
\begin{lemma}\label{Mon}
Assume $\frac{n}{n-4} < p < \frac{n+4}{n-4}$  and that $w$ is a nonnegative  $C^4$ solution of \eqref{W01}.  Then, $E(r; w)$ is non-increasing in $t \in (-\infty, 0)$. Furthermore, we have 
\begin{equation}\label{Mon01}
\frac{d}{dt} E(t; w)= K_3\int_{\mathbb{S}^{n-1}} \left[ \left( \partial_t^2 w \right)^2 + |\partial_t \nabla_\theta w|^2 \right] - K_1 \int_{\mathbb{S}^{n-1}} \left( \partial_t w \right)^2.  
\end{equation}
\end{lemma}

 \begin{remark}\label{RMY0}
 An analogous monotonicity formula has been derived by the author and Zou \cite{YZ} to study  isolated singularities for a fractional equation.  Ghergu-Kim-Shahgholian \cite{GKS} also obtained a similar monotonicity formula for a  second order semilinear elliptic system with power-law nonlinearity.     
 \end{remark}

\bp
Note that 
$$
\aligned
\partial_t^{(4)}w \partial_t w & = \partial_t \left( \partial_t^{(3)}w \partial_t w \right) -\partial_t^{(3)}w \partial_t^2 w \\
& = \partial_t \Big( \partial_t^{(3)}w \partial_t w - \frac{1}{2} \left( \partial_t^2 w \right)^2 \Big),  \\
\partial_t^{(3)}w \partial_t w & = \partial_t \left( \partial_t^2 w \partial_t w \right) - \left( \partial_t^2 w \right)^2, \\
\partial_t^2w  \partial_t w & = \frac{1}{2} \partial_t \left( \partial_t w \right)^2. 
\endaligned
$$
Therefore, multiplying Eq. \eqref{W01} by $\partial_t w$ and integrating by parts on $\mathbb{S}^{n-1}$, we get  
$$
\frac{d}{dt} E(t; w)= K_3\int_{\mathbb{S}^{n-1}}  \left[ \left( \partial_t^2 w \right)^2 + |\partial_t \nabla_\theta w|^2 \right]   - K_1 \int_{\mathbb{S}^{n-1}} \left( \partial_t w \right)^2. 
$$
By Lemma \ref{L-S301}, we have $K_1 >0$ and $K_3 <0$. Hence $\frac{d}{dt} E(t; w) \leq 0$ and we finish the proof.  
\ep

\begin{lemma}\label{Bou}
Let $w$ be a nonnegative $C^4$ solution of \eqref{W01} with $1 < p < \frac{n+4}{n-4}$. Then $w$, $\partial_t w$, $\partial_t^2 w$, $\partial_t^{(3)} w$, $\Delta_\theta w$ and $|\nabla_\theta w|$ are uniformly bounded in $(-\infty, -\ln 2) \times \mathbb{S}^{n-1}$.  
\end{lemma}
\bp Define
$$
u(x) = |x|^{-\frac{4}{p-1}} w(t, \theta),
$$
where $t=\ln |x|$ and $\theta =\frac{x}{|x|}$.  Then $u$ is a nonnegative solution of \eqref{Bi}. By Lemma \ref{Le201}, we know that $w$ is uniformly bounded. By Corollary  \ref{Co201} we have
$$
\aligned
| \partial_t w| + |\nabla_\theta w|  & \leq  C \sum_{i=0}^1 |x|^{\frac{4}{p-1} + i} |\nabla_x^i u| \leq C, \\
| \partial_t^2  w| + |\Delta_\theta w|  & \leq  C \sum_{i=0}^2 |x|^{\frac{4}{p-1} + i} |\nabla_x^i u| \leq C, \\
|\partial_t^{(3)} w|  & \leq C \sum_{i=0}^3 |x|^{\frac{4}{p-1} + i} |\nabla_x^i u| \leq C. 
\endaligned
$$
Thus the desired conclusion follows. 
\ep

Assume $\frac{n}{n-4} < p < \frac{n+4}{n-4}$, from Lemmas \ref{Mon} and \ref{Bou} we deduce that
the limit $ \lim_{t \to -\infty} E(t; w)$
exists.  Let $u$ be a nonnegative solution of \eqref{Bi}, we define
\begin{equation}\label{EE}
\widetilde{E}(r; u):=E(t; w), 
\end{equation}
where $t=\ln r$ and $w$ is defined as in \eqref{Set2}. Then we have 
$$
\widetilde{E}(0; u):=\lim_{r\to 0^+} \widetilde{E}(r; u)=\lim_{t \to -\infty} E(t; w). 
$$
For any $\lambda >0$, define
$$
u^\lambda(x):=\lambda^{\frac{4}{p-1}} u(\lambda x). 
$$
Then $u^\lambda$ is also a nonnegative solution of \eqref{Bi} in $B_{1/\lambda}  \backslash \{0\}$. 
Moreover, we have 
$$
\aligned
\widetilde{E}(r; u^\lambda) & = E(t; w(\cdot + \ln \lambda, \cdot))\\
&= E(t+ \ln \lambda, w)\\
&= E(\lambda r, u). 
\endaligned
$$
That is, we get the following  scaling invariance
\begin{equation}\label{Sac}
\widetilde{E}(r; u^\lambda)= \widetilde{E}(\lambda r; u). 
\end{equation} 
\begin{lemma}\label{Lim}
Let $u \in C^4(B_1 \backslash \{0\})$  be a nonnegative solution of \eqref{Bi} with $\frac{n}{n-4} < p <\frac{n+4}{n-4}$. Then either 
$$\lim_{|x| \to 0} |x|^{\frac{4}{p-1}} u(x)=0$$ 
or 
$$\lim_{|x| \to 0} |x|^{\frac{4}{p-1}} u(x)=K_0^{\frac{1}{p-1}},$$ 
where $K_0$ is given by \eqref{K0}. 
\end{lemma}
\bp
First we compute the possible values of $\widetilde{E}(0; u)$. By Lemma \ref{Le201}, $u^\lambda$ are uniformly bounded in $C^{4, \alpha} (K)$ on every compact set $K \subset B_{1 / 2\lambda}  \backslash \{0\}$, with some $0< \alpha <1$.  Therefore,  there exists a  nonnegative function $u^0 \in C^4(\mathbb{R}^n \backslash \{0\})$, such that up to a subsequence of $\lambda \to 0$, $u^\lambda$ converges to $u^0$ in $C_{loc}^4(\mathbb{R}^n \backslash \{0\})$. Further, $u^0$ satisfies
$$
\Delta^2 u =u^p ~~~~~~~~~\textmd{in} ~  \mathbb{R}^n \backslash \{0\}.
$$
By Lemma \ref{Non589}, we have $-\Delta u^0 \geq 0$ in $\mathbb{R}^n \backslash \{0\}$. The maximum principle gives that either 
$$
u^0 \equiv 0 ~~~~~~~~~~  \textmd{in}  ~ \mathbb{R}^n \backslash \{0\}
$$
or 
$$
u^0 >0 ~~~~~~~~~~  \textmd{in}  ~ \mathbb{R}^n \backslash \{0\}.
$$
Therefore, by Theorem 4.2 in \cite{Lin}, $u^0$ is radially symmetric with respect to the origin 0. 
Moreover,  by the scaling invariance of $\widetilde{E}$, we have for any $r >0$ that 
\begin{equation}\label{Con}
\widetilde{E}(r; u^0)= \lim_{\lambda \to 0} \widetilde{E}(r; u^\lambda)=\lim_{\lambda \to 0} \widetilde{E}(\lambda r; u) =\widetilde{E}(0; u). 
\end{equation}
Let
$$
w^0(t) := |x|^{\frac{4}{p-1}} u^0(|x|), ~~~~~t=\ln |x|. 
$$
Then $w^0$ satisfies 
\begin{equation}\label{W023}
\frac{d^4}{dt^4}w^0  + K_3 \frac{d^3}{dt^3} w^0 + K_2 \frac{d^2}{dt^2} w^0 + K_1 \frac{d}{dt}  w^0 + K_0 w^0 = (w^0)^p  ~~~~ \textmd{in} ~  \mathbb{R}.  
\end{equation}
From \eqref{Con}, $E(t; w^0)=\widetilde{E}(r; u^0)$  is a constant. By Lemma \ref{Mon}, 
$$
\frac{d}{dt} E(t; w^0)= |\mathbb{S}^{n-1}| \left[ K_3 \left( \frac{d^2}{dt^2} w^0 \right)^2 - K_1  \left( \frac{d}{dt} w^0 \right)^2 \right] \equiv 0. 
$$
Since $K_3 <0$ and $K_1 >0$, we get  that $\frac{d}{dt} w^0 \equiv 0$ in $\mathbb{R}$ and then $w^0$ is a constant.  By \eqref{W023}, either $w^0=0$ or $w^0=K_0^{\frac{1}{p-1}}$. Hence, by \eqref{Con} we obtain 
$$
\widetilde{E}(0; u)  \in \left\{ 0, \left(\frac{1}{2} - \frac{1}{p+1}\right) K_0^{\frac{p+1}{p-1}} |\mathbb{S}^{n-1}| \right\}.
$$

If $\widetilde{E}(0; u) =0$, then $u^0 \equiv 0$. Since this function $u^0$ is unique, we conclude that $u^\lambda \to 0$ for any sequence of $\lambda \to 0$, in $C^{4} (K)$ on every compact set $K \subset \mathbb{R}^n  \backslash \{0\}$.  Therefore, we easily get
$$
\lim_{|x| \to 0} |x|^{\frac{4}{p-1}} u(x) =0. 
$$

\vskip0.10in

If $\widetilde{E}(0; u) =\left(\frac{1}{2} - \frac{1}{p+1}\right) K_0^{\frac{p+1}{p-1}} |\mathbb{S}^{n-1}|$, then we have 
$$
u^0(x) \equiv  K_0^{\frac{1}{p-1}} |x|^{-\frac{4}{p-1}}.
$$
In this case the function $u^0$ is also unique, so we obtain that $u^\lambda \to K_0^{\frac{1}{p-1}} |x|^{-\frac{4}{p-1}}$ for any sequence of $\lambda \to 0$, in $C^{4} (K)$ on every compact set $K \subset \mathbb{R}^n  \backslash \{0\}$.  In particular, we have
$$
u^\lambda(x) \to K_0^{\frac{1}{p-1}} ~~~~~~ \textmd{as} ~ \lambda \to 0
$$
in $C (\mathbb{S}^{n-1})$.  We quickly get
$$
\lim_{|x| \to 0} |x|^{\frac{4}{p-1}} u(x) = K_0^{\frac{1}{p-1}}.
$$
This completes the proof. 
\ep

\begin{lemma}\label{Rem987}
Assume $\frac{n}{n-4} < p < \frac{n+4}{n-4}$ and that $u \in C^4(B_1 \backslash \{0\})$ is a nonnegative solution of \eqref{Bi}. If 
\begin{equation}\label{Assu90}
\lim_{|x| \to 0} |x|^{\frac{4}{p-1}} u(x)=0,
\end{equation}
then
\begin{equation}\label{Re789}
\int_{ \{ |x|\leq 1/2 \} } u^{(p-1)n/4} < +\infty. 
\end{equation}
\end{lemma}
\bp
Set
$$
\varphi(x) =|x|^{-\frac{n-4}{p-1}\left( p- \frac{n}{n-4}\right)}.
$$
We recall that in radial coordinates $r=|x|$, we have
\begin{equation}\label{Ra-C0}
\aligned
\Delta^2 \varphi(r) = & \varphi^{(4)}(r) + \frac{2(n-1)}{r} \varphi^{(3)}(r) + \frac{(n-1)(n-3)}{r^2} \varphi^{\prime\prime}(r) \\
& - \frac{(n-1)(n-3)}{r^3} \varphi^{\prime}(r)
\endaligned
\end{equation}
Denote
$$
\gamma := -\frac{n-4}{p-1}\left( p- \frac{n}{n-4}\right) < 0. 
$$
Direct calculations show that 
$$
\aligned
\Delta^2 \varphi  & =  \Big[ \gamma(\gamma-1)(\gamma-2)(\gamma-3) + 2\gamma (n-1) (\gamma-1)(\gamma-2) \\
& ~~~~~ + \gamma(n-1)(n-3)(\gamma-1) -\gamma(n-1)(n-3) \Big] r^{\gamma-4} \\
& = \gamma(\gamma-2)\Big[ (\gamma-1)(\gamma-3) + 2(n-1)(\gamma-1) +(n-1)(n-3) \Big] r^{\gamma-4} \\
& = \gamma(\gamma-2)\Big[ (\gamma-1) (\gamma+n-4) + (n-1)(\gamma+n-4) \Big] r^{\gamma-4}\\
& = \gamma(\gamma-2)(\gamma + n-2)(\gamma+n-4) r^{\gamma-4}. 
\endaligned
$$
Since $\gamma +n -4=\frac{4}{p-1} >0$, we have  
$$
A:= \gamma(\gamma-2)(\gamma + n-2)(\gamma+n-4) >0. 
$$
That is, we obtain
\begin{equation}\label{AA}
\frac{\Delta^2 \varphi}{\varphi} =\frac{A}{|x|^4} ~~~~~~~ \textmd{in} ~ \mathbb{R}^n \backslash  \{0\}. 
\end{equation}
For small $\epsilon>0$, let $\zeta_\epsilon$ be a smooth cut-off function satisfying 
\begin{equation}\label{kesi}
\zeta_\epsilon(x)=
\begin{cases}
1 ~~~~~~~ & \textmd{for} ~  \epsilon \leq |x| \leq \frac{1}{2}, \\
0 ~~~~~~~ & \textmd{for} ~ |x|\leq \frac{\epsilon}{2}, ~ |x|\geq \frac{3}{4}
\end{cases}
\end{equation}
and 
\begin{equation}\label{kesi01}
|\nabla^k \zeta_\epsilon (x)| \leq C \epsilon^{-k} ~~~~~~~  \textmd{for} ~ k=1, 2, 3, 4. 
\end{equation}
Using $\zeta_\epsilon \varphi$ as a test function in \eqref{Bi} and integrating by parts we obtain 
\begin{equation}\label{Int-09}
\int_{B_1} \zeta_\epsilon u \varphi \left( \frac{\Delta^2 \varphi}{\varphi} -u^{p-1} \right) = - \int_{B_1} u F(\zeta_\epsilon, \varphi), 
\end{equation}
where
$$
\aligned
F(\zeta_\epsilon, \varphi)= & 4 \nabla \zeta_\epsilon \cdot \nabla \Delta\varphi + 2\Delta \zeta_\epsilon \Delta\varphi +4 \sum_{i,j=1}^n (\zeta_\epsilon)_{x_ix_j} \varphi_{x_i x_i} \\
& + 4 \nabla \Delta\zeta_\epsilon \cdot \nabla \varphi + \varphi \Delta^2 \zeta_\epsilon. \\
\endaligned
$$
By \eqref{kesi}, \eqref{kesi01} and Lemma \ref{Le201}, we estimate 
$$
\aligned
\left| \int_{B_1} u F(\zeta_\epsilon, \varphi) \right|  & \leq \left| \int_{ \{\frac{1}{2} \leq |x| \leq \frac{3}{4} \} } u F(\zeta_\epsilon, \varphi) \right| + \left| \int_{  \{ \frac{\epsilon}{2} \leq |x| \leq \epsilon \} } uF(\zeta_\epsilon, \varphi) \right| \\
& \leq C_1 + C_2 \int_{  \{ \frac{\epsilon}{2} \leq |x| \leq \epsilon \} } u \bigg[ \frac{1}{\epsilon} |x|^{\gamma-3} + \frac{1}{\epsilon^2} |x|^{\gamma-2} \\
& ~~~~~  + \frac{1}{\epsilon^3} |x|^{\gamma-1} + \frac{1}{\epsilon^4}|x|^{\gamma}  \bigg] \\
& \leq C_1 + C_2 \epsilon^n \epsilon^{\gamma-4} \epsilon^{ - \frac{4}{p-1}} \leq C_1 + C_2 < +\infty,  \\
\endaligned
$$
where $C_1=C_1(p, n, u)$ and $C_2=C_2(p, n)$ are two positive constants (independent of $\epsilon$).  Hence
\begin{equation}\label{Int-765}
\int_{B_1} \zeta_\epsilon u \varphi \left( \frac{\Delta^2 \varphi}{\varphi} -u^{p-1} \right) \leq C_2 < +\infty 
\end{equation}
uniformly in $\epsilon$. By the assumption \eqref{Assu90}, 
$$ 
u^{p-1}(x)=o(1) |x|^{-4} ~~~~~~~~ \textmd{as} ~ |x| \to 0.
$$
This together with \eqref{AA} and \eqref{Int-765} gives 
$$
\int_{B_1} \zeta_\epsilon u |x|^{\gamma-4} \leq C_3 < +\infty, 
$$
where $C_3$ is a positive constant  independent of $\epsilon$. Therefore, by Lemma \ref{Le201}, 
$$
\aligned
\int_{ \{ \epsilon \leq |x| \leq \frac{1}{2} \} } u^{(p-1)n / 4} &=  \int_{ \{ \epsilon \leq |x| \leq \frac{1}{2} \} } u u^{((p-1)n -4) / 4 } \\
& \leq C(p, n) \int_{ \{ \epsilon \leq |x| \leq \frac{1}{2} \} } u |x|^{- ( (p-1)n -4 ) / (p-1)} \\
& =C(p, n) \int_{ \{ \epsilon \leq |x| \leq \frac{1}{2} \} } u |x|^{\gamma -4} \\ 
& \leq  C(p, n) \int_{B_1} \zeta_\epsilon u |x|^{\gamma-4} \leq C(p, n)  C_3 < +\infty. 
\endaligned
$$ 
Letting $\epsilon \to 0$, we get \eqref{Re789} by the dominated convergence theorem.
\ep

Now we give a new method to obtain the removable  singularity theorem. For our fourth order equation \eqref{Bi}, the classical methods based on the maximum principle to  second order problems (such as see \cite{A,CJSX,CL2,KMPS}) fail.  We remark that our method also apply to higher order equations. This method is based on the following Regularity Lifting Theorem from Chen-Li \cite{CL-Book}. 

\vskip0.10in

Let $V$ be a Hausdorff topological vector space. Suppose there are two extended norms (i.e., the norm of an element in $V$ might be infinity) defined on $V$, 
$$
\|\cdot\|_X, ~ \|\cdot\|_Y : V \to [0, +\infty]. 
$$
Let 
$$
X:=\{v \in V : \|v\|_X < +\infty\} ~~~~ \textmd{and} ~~~~ Y:=\{v \in V : \|v\|_Y < +\infty\}. 
$$
Assume that spaces $X$ and $Y$ are complete under the corresponding norms and the convergence in $X$ or in $Y$ implies the convergence in $V$.
\begin{theorem}\label{CLBk} 
(\cite{CL-Book}, Theorem 3.3.1) Let $T$ be a contraction map from $X$ into itself and from $Y$ into itself. Assume that $f\in X$ and that there exists a function $g\in Z:=X\cap Y$ such that $f=Tf +g$ in $X$. Then $f$ also belongs to $Z$. 
\end{theorem}

\begin{remark}
We usually choose $V$ to be the space of distributions, and $X$ and $Y$ to be function spaces, for instance, $X=L^p$ and $Y=W^{1, q}$. 
\end{remark}

 Next we will use this Regularity Lifting Theorem to prove a removable singularity result.  
\begin{lemma}\label{Rem345}
Assume $\frac{n}{n-4} < p < \frac{n+4}{n-4}$ and that $u \in C^4(B_1 \backslash \{0\})$  is a nonnegative solution of \eqref{Bi}. If 
\begin{equation}\label{Rem-093}
\int_{ \{ |x|\leq 1/2 \} } u^{(p-1)n/4  } < +\infty,
\end{equation}
then the singularity at $x=0$ is removable, i.e., $u(x)$ can be extended to a $C^4$ solution of \eqref{Bi} in the entire ball $B_1$. 
\end{lemma}
\begin{remark}
For the second order equation \eqref{Sec}, a similar result for removable singularity was proved by Gidas-Spruck \cite{GS}. However, their proof is based on a double application of the De Giorgi-Nash-Moser bootstrap arguments, which cannot be applied to our fourth order problem \eqref{Bi}.  
\end{remark}
\bp
Let $G_2(x, y)$ be the Green's function of $\Delta^2$ in $B_{1/2}$ with homogeneous Dirichlet boundary conditions, Then, for each fixed $y\in B_1$, $G_2(\cdot, y)$ is a distributional solution of
$$
\begin{cases}
\Delta^2 G_2(\cdot, y) =\delta(\cdot - y) ~~~~~~~~ & \textmd{in} ~ B_{1/2},\\
G_2(\cdot, y)=\frac{\partial G_2(\cdot, y)}{\partial \nu} =0 ~~~~~~~& \textmd{on} ~ \partial B_{1/2},
\end{cases}
$$
and there exists positive constant $C_n$ such that
$$
0< G_2(x, y) \leq \Gamma_2(|x - y|):= C_n |x-y|^{4-n}~~~~~~~ \textmd{for}  ~x, y \in B_{1/2},  x \neq y. 
$$
Define 
$$
v(x):= -u(x)+ \int_{B_{1/2}} G_2(x, y)u^p(y)dy, ~~~~~ x \in B_{1/2}, 
$$
then $v \in L^1(B_{1/2})$. Moreover, by Lemma \ref{Weak},  $v$ satisfies 
$$
\Delta^2 v =0 ~~~~~~~ \textmd{in} ~ B_{1/2}
$$
in the distributional sense.  Using Theorem 7.23 in \cite{Mi}, we get $v \in L_{loc}^\infty (B_{1/2})$. 

\vskip0.10in

Now we split the right hand side of \eqref{Bi} into two parts:
$$
u^p=u^{p-1} u:=a(x) u.
$$
Then, by the assumption \eqref{Rem-093}, $a(x) \in L^\frac{n}{4}(B_{1/2})$. For any positive number $L>0$, let
$$
a_L(x)=
\begin{cases}
a(x)~~~~~~~~ & \textmd{if} ~ |a(x)| \geq L, \\
0~~~~~~~ & \textmd{otherwise},
\end{cases}
$$
and
$$
a_M(x)=a(x) - a_L(x).
$$
Define the linear operator
$$
(T_L w)(x)=\int_{B_{1/4}} G_2(x, y) a_L(y) w(y) dy.
$$
Then $u$ satisfies the equation
\begin{equation}\label{Inte953}
u(x)= (T_L u)(x) + F_L(x) ~~~~~~~~ \textmd{in}  ~ B_{1/4}, 
\end{equation}
where
$$
F_L(x)=\int_{B_{1/4}}G_2(x, y)a_M(y)u(y) dy - v(x) +h(x) 
$$
and
$$
h(x)=\int_{ \{\frac{1}{4} \leq |y| \leq \frac{1}{2} \} } G_2(x, y) u^p(y)dy.
$$
Note that
$$
\aligned
|h(x)| & \leq C\int_{ \{\frac{1}{4} \leq |y| \leq \frac{1}{2} \} } G_2(x, y)dy \leq C \int_{ \{\frac{1}{4} \leq |y| \leq \frac{1}{2} \} } |x -y|^{4-n} dy\\
&  \leq C\int_{B_1} |y|^{4-n} dy \leq C~~~~~~~~ \textmd{for} ~  \textmd{all}  ~ x \in ~ B_{1/4}. 
\endaligned
$$
 Hence $v, h \in L^\infty(B_{1/4})$. 

\vskip0.10in

We will prove that, for any $\frac{n}{n-4} < q < \infty$,

\vskip0.10in

(1) $T_L$ is a contracting operator from $L^q(B_{1/4})$ to $L^q(B_{1/4})$ for $L$ large. 

\vskip0.10in

(2) $F_L \in L^q(B_{1/4})$. 

\vskip0.10in

Then, by the Regularity Lifting Theorem \ref{CLBk}, we obtain $u\in L^q(B_{1/4})$ for any $\frac{n}{n-4} < q < \infty$. 
 
\vskip0.10in
 
(1) {\it The estimate of the operator $T_L$.}
 
\vskip0.10in

For any $\frac{n}{n-4} < q <  \infty$, there exists $1 < r < \frac{n}{4}$ such that
$$
\frac{1}{q}=\frac{1}{r} -\frac{4}{n}. 
$$ 
By Hardy-Littlewood-Sobolev inequality and H\"older inequality, we have
$$
\aligned
 \|T_L w\|_{L^q(B_{1/4})} & \leq \|\Gamma_2 \ast a_L w\|_{L^q(\mathbb{R}^n)} \leq C \|a_L w\|_{L^r({B_{1/4}})} \\
 &  \leq \|a_L\|_{L^{\frac{n}{4}} (B_{1/4})}  \|w\|_{L^q(B_{1/4})}. 
 \endaligned
$$
Since $a(x) \in L^{\frac{n}{4}}(B_{1/4})$, we can choose $L$ sufficiently large, such that
$$
\|a_L\|_{L^{\frac{n}{4}} (B_{1/4})} \leq \frac{1}{2}. 
$$
Therefore,  $T_L : L^q(B_{1/4}) \to L^q(B_{1/4})$ is a contracting operator for $L$ large. 
\vskip0.10in
 
(2) {\it The integrability of the function $F_L(x)$.}
 
\vskip0.10in

Obviously, we only need to show that, for any $\frac{n}{n-4} < q < \infty$,
$$
F_L^1(x) := \int_{B_{1/4}}G_2(x, y)a_M(y)u(y) dy \in L^q(B_{1/4}). 
$$
Since $a_M(x)$ is a bounded function, we have
$$
\|F_L^1\|_{L^q (B_{1/4})} \leq \|a_M u\|_{L^r(B_{1/4})} \leq C\|u\|_{L^r(B_{1/4})}.
$$
By the assumption \eqref{Rem-093}, $u \in L^r(B_{1/4})$ for any $1 < r \leq \frac{(p-1)n}{4}$. Note that 
$$
q=\frac{(p-1) n}{ 4(2-p)}~~~~ \textmd{if} ~ r=\frac{(p-1)n}{4}. 
$$
Hence,  we conclude that, for the following values of $q$,                  
$$
\begin{cases}
1< q < \infty ~~~~~~ & \textmd{if} ~ p\geq 2, \\
1 < q \leq \frac{(p-1)n}{4(2-p)} & \textmd{if} ~ p <  2, 
\end{cases}
$$
$F_L(x) \in L^q(B_{1/4})$. 

\vskip0.10in

Using the Regularity Lifting Theorem \ref{CLBk}, we obtain
$$
\begin{cases}
u \in L^q(B_{1/4})~~~ & \textmd{for} ~ \textmd{any} ~1< q < \infty ~~~~~ \textmd{if} ~ p\geq 2, \\
u \in L^q(B_{1/4})~~~ & \textmd{for} ~ \textmd{any} ~1 < q \leq \frac{(p-1)n}{4(2-p)} ~~~~ \textmd{if} ~ p <  2. 
\end{cases}
$$
Now we note that  from the starting point where $u \in L^\frac{(p-1)n}{4}(B_{1/4})$, we get 
$$
u \in L^r(B_{1/4})~~~~~~  \textmd{with}~ r = \frac{(p-1)n}{4(2-p)},~~ p<2. 
$$
By a similar argument as above, we  get
$$
\begin{cases}
u \in L^q(B_{1/4})~~~ & \textmd{for} ~ \textmd{any} ~1< q < \infty ~~~~~ \textmd{if} ~ p\geq \frac{3}{2}, \\
u \in L^q(B_{1/4})~~~ & \textmd{for} ~ \textmd{any} ~1 < q \leq \frac{(p-1)n}{4(3-2p)} ~~~~ \textmd{if} ~ p <  \frac{3}{2}. 
\end{cases}
$$
Hence by iteration we have for $k=1, 2, \cdots$,  
$$
\begin{cases}
u \in L^q(B_{1/4})~~~ & \textmd{for} ~ \textmd{any} ~1< q < \infty ~~~~~ \textmd{if} ~ p\geq \frac{k+1}{k}, \\
u \in L^q(B_{1/4})~~~ & \textmd{for} ~ \textmd{any} ~1 < q \leq \frac{(p-1)n}{4[(k+1) - kp]} ~~~~ \textmd{if} ~ p <  \frac{k+1}{k}. 
\end{cases}
$$  
This implies that for any fixed dimension $n$, a finite number of iterations gives
$$
u \in L^q(B_{1/4}) ~~~~~  \textmd{for} ~ \textmd{any} ~1< q < \infty. 
$$

\vskip0.10in
 Finally, By H\"older inequality, we have
 $$
 \int_{B_{1/4}} G_2(x, y) u^p(y)dy \in L^\infty(B_{1/4}), 
 $$
From this and \eqref{Inte953} we easily deduce   that $u \in L^\infty(B_{1/4})$. By estimates of elliptic equations, $u(x)$ is smooth at 0. Therefore 0 is a removable singularity.  
\ep

\vskip0.20in

\noindent{\it Proof of Theorem \ref{T-Bi}.}  The proof of Theorem \ref{T-Bi} is now just a combination of Lemmas \ref{Lim}, \ref{Rem987} and \ref{Rem345}.
\hfill$\square$

\vskip0.20in

\noindent{\bf Acknowledgments.} The author would like to thank Professor Sun-Yung A. Chang for many helpful discussions and comments.  The author would also like to thank his advisor Professor Wenming Zou for his constant support and encouragement. 
This work was done during the author's visit to Princeton University.  He thanks Tsinghua University for funding his visit and thanks the Department of Mathematics at Princeton University for kind hospitality.

\end{document}